\documentclass[11pt]{amsart}
\usepackage{mesfonts}

\def\S{\mathcal S}

\def\p{\mathbb{P}}

\def\penn{\pen}

\def\crit{{\mathrm{Crit}}}

\def\S{\mathcal{S}}

\newcommand{\kf}[2]{\|#1\|_{(2,#2)}}
\newcommand{\KF}{\mathrm{K\!F}}
\newcommand{\RSC}{\mathrm{RSC}}
\newcommand{\RSCI}{\mathrm{RSCI}}
 
\def\penn{\overline{\pen}}

\RequirePackage{amssymb}
\newtheorem{lemma}{Lemma}
\newtheorem{theorem}{Theorem}
\newtheorem{proposition}{Proposition}
\newtheorem{corollary}{Corollary}

\newtheorem{fact}{Fact}

\usepackage{graphicx}
\RequirePackage{amssymb}

\begin{document}

\title{Low rank multivariate regression}
\author{Christophe Giraud}
\address{CMAP, UMR CNRS 7641, Ecole Polytechnique, Route de Saclay, 91128 Palaiseau Cedex, France} 
\email{christophe.giraud@polytechnique.edu}
\date{September 2010, revision April 2011}
\begin{abstract}We consider in this paper the multivariate regression problem, when the target regression matrix $A$ is close to a low rank matrix. Our primary interest is in on the practical  case where the variance of the noise is unknown.
Our main contribution is to propose in this setting a criterion to select among a family of low rank estimators and prove a non-asymptotic oracle inequality for the resulting estimator. We also investigate the easier case where the variance of the noise is known and outline that the penalties appearing in our criterions are minimal (in some sense). These penalties involve the expected value of Ky-Fan norms of some random matrices. These quantities can be evaluated easily in practice and upper-bounds can be derived from recent results in random matrix theory.
\end{abstract}

\maketitle

\section{Introduction}
We build on ideas introduced in a recent paper of Bunea, She and Wegkamp~\cite{BuneaSheWegkamp,Oct2010}
for the multivariate regression problem
\begin{equation}\label{model}
Y=XA+\sigma E
\end{equation}
where $Y$ is a $m\times n$ matrix of response variables, $X$ is a $m\times p$ matrix of predictors, $A$ is $p\times n$ matrix of regression coefficients and $E$ is a $m\times n$ random matrix with i.i.d. entries.
We assume for simplicity that the entries $E_{i,j}$ are standard Gaussian,
yet  all the results can be extended to the case where the entries are sub-Gaussian.

An important issue in multivariate regression is to estimate $A$ or $XA$ when the matrix $A$ has a low rank or can be well approximated by a low rank matrix, see Izenman~\cite{Izenman08}. In this case, a small number of linear combinations of the predictors catches most of the non-random variation of the response $Y$. 
This framework arises in many applications, among which analysis of fMRI image data~\cite{Harrisonetal}, analysis of EEG data decoding~\cite{Andersonetal}, neural response modeling~\cite{Brownetal} or genomic data analysis~\cite{BuneaSheWegkamp}. 

When the variance $\sigma^2$ is known, the strategy developed by Bunea {\it et al.}~\cite{BuneaSheWegkamp} for estimating $A$ or $XA$  is the following.
Writing   $\|.\|$ for the Hilbert-Schmidt norm and $\widehat A_{r}$ for the minimizer of 
$\|Y-X\widehat A\|$
over the matrices $\widehat A$ of rank at most $r$, the matrix $XA$ is estimated by $X\widehat A_{\hat r}$, where $\hat r$ minimizes the criterion
\begin{equation}\label{critvarianceconnue}
\crit_{\sigma^2}(r)=\|Y-X\widehat A_{r}\|^2+\pen_{\sigma^2}(r)\sigma^2.
\end{equation}
Bunea {\it et al.}~\cite{BuneaSheWegkamp} considers a penalty $\pen_{\sigma^2}(r)$ linear in $r$ and provides clean non-asymptotic bounds on $\|X\widehat A_{\hat r}-XA\|^2$, on $\|\widehat A_{\hat r}-A\|^2$ and on the probability that the estimated rank $\hat r$ coincides with the rank of $A$.

Our main contribution is to propose and analyze a criterion to handle the case where $\sigma^2$ is unknown.
Our theory requires no assumption on the design matrix $X$ and applies in particular when the sample size $m$ is smaller than the number of covariates $p$.
We also exhibit a minimal sublinear penalty for the Criterion~(\ref{critvarianceconnue}) for the case of known variance.

Let us denote by $q$ the rank of $X$ and by $G_{q\times n}$ a $q\times n$ random matrix with i.i.d. standard Gaussian entries.
The penalties that we introduce  involve the expected value of the Ky-Fan $(2,r)$-norm of the random matrix $G_{q\times n}$, namely
$$\S_{q\times n}(r)=\E\cro{\kf{G_{q\times n}}{r}},\quad\textrm{where}\;\      \kf{G_{q\times n}}{r}^2=\sum_{k=1}^r\sigma^2_{k}(G_{q\times n})$$
and where $\sigma_{k}(G_{q\times n})$ stands for the $k$-th largest singular value of $G_{q\times n}$. The term $\S_{q\times n}(r)$ can be evaluated by Monte Carlo and for $q,n$ large enough an accurate approximation of $\S_{q\times n}(r)$  is derived from the Marchenko-Pastur distribution, see Section~2.

 For the case of unknown variance,
we prove a non-asymptotic oracle-like inequality for the criterion 
\begin{equation}\label{critunknown2}
\crit(r)=\log(\|Y-X\hat A_{r}\|^2)+\pen(r).
\end{equation}
with
$$\pen(r)\geq -\log\pa{1-K\,{\S_{q\times n}(r)^2\over nm-1}},\quad\textrm{with}\ K>1.$$
The latter constraint on the penalty is shown to be minimal (in some sense).
In addition, we also consider the case where $\sigma^2$ is known and show that the penalty 
$\pen(r)=\S_{q\times n}(r)^2$
is minimal  for the Criterion~(\ref{critvarianceconnue}).

The study of multivariate regression with rank constraints dates back to Anderson~\cite{Anderson} and Izenman~\cite{Izenman}. The question of rank selection has only been recently addressed by Anderson~\cite{Anderson} in an asymptotic setting (with $p$ fixed) and by Bunea {\it et al.}~\cite{BuneaSheWegkamp, Oct2010} 
in an non-asymptotic framework. We refer to the latter article for additional references. In parallel, a series of recent papers 
study the estimator $\widehat A_{\lambda}^{\ell^1}$ obtained by minimizing
$$\|Y-X\widehat A\|^2+\lambda  \sum_{k} \sigma_{k}(\widehat A)$$
see among others Yuan {\it et al.}~\cite{Yuanetal}, Bach~\cite{Bach}, Neghaban and Wainwright~\cite{NeghabanWainwright}, 
Lu {\it et al.}~\cite{Luetal}, Rohde and Tsybakov~\cite{RohdeTsybakov}, Koltchinskii {\it et al.}~\cite{Koltchinskii}. Due to the "$\ell^1$" penalty 
$\sum_{k} \sigma_{k}(\widehat A)$, the estimator $\widehat A_{\lambda}^{\ell^1}$ has a small rank for $\lambda$ large enough and it is proven to have good statistical properties under some hypotheses on the design matrix $X$. 
We refer to Bunea {\it et al.}~\cite{Oct2010} for a detailed analysis of the similarities and the differences  between  $\widehat A_{\lambda}^{\ell^1}$ and their estimator.

Our paper is organized as follows. In the next section, we give a few results on $\S_{q\times n}(r)$ and on the estimator $X\widehat A_{r}$. In Section~3, we analyze the case where the variance $\sigma^2$ is known, which gives us a benchmark for the Section~4 where the case of unknown variance is tackled. In Section~5, we comment on the extension of the results to the case of sub-Gaussian errors and 
 we outline that our theory provides a theoretically grounded criterion (in a  non-asymptotic framework)
to select the number $r$ of components to be kept in a principal component analysis. Finally, we carry out  an empirical study in Section~6 and prove the main results in Section~7.

\subsection*{{\tt R}-code} The estimation procedure described in sections~4 and~7 has been implemented in {\tt R}. We provide the {\tt R}-code (with a short notice) at the following URL : \\*
{\tt http://www.cmap.polytechnique.fr/$\sim$giraud/software/KF.zip}

\subsection*{What is new here?}
The primary purpose of the first draft of the present paper~\cite{Giraud} was to provide complements to the paper of Bunea {\it et al.}~\cite{BuneaSheWegkamp} in the two following directions:
\begin{itemize}
\item to propose a selection criterion for the case of unknown variance,
\item to give some tighter results for Gaussian errors.
\end{itemize}
During the reviewing process of the first draft of this paper,  Bunea, She and Wegkamp wrote an augmented version of their paper~\cite{Oct2010} were they also investigate these two points. Let us comment briefly on the overlap between the results of these two simultaneous works~\cite{Giraud,Oct2010}. Let us start with the main contribution of our paper, which is to provide a selection criterion for the case of unknown variance. In Section~2.4 of~\cite{Oct2010}, the authors propose and analyze a criterion to handle the case of unknown variance in the setting where the rank $q$ of $X$ is {\it strictly} smaller than the sample size $m$. In this favorable case, the variance $\sigma^2$ can be conveniently estimated by
$$\hat \sigma^2={\|Y-PY\|^2\over mn-qn},\quad\textrm{with } P\ \textrm{the orthogonal projector onto the range of }X,$$
which has the nice feature to be an unbiased estimator of $\sigma^2$  independent of the collection of estimators $\{\widehat A_{r},\ r=0,\ldots,q\}$. 
Plugging this estimator $\hat \sigma^2$ in the Criterion~(\ref{critvarianceconnue}), Bunea {\it et al.}~\cite{Oct2010} proves a nice oracle bound. This approach no more applies in the general case where the rank of $X$ can be as large as $m$, which is very likely to happen when the number  $p$  of covariates is larger than the sample size~$m$. We provide in Section~4 an oracle inequality for the Criterion~(\ref{critunknown2}) with no restriction on the rank of $X$.
 
 Concerning the case of known variance : the final paper of Bunea  {\it et al.}~\cite{Oct2010} proposes for Gaussian errors the penalty $\pen_{\sigma^2}(r)=Kr(\sqrt{q}+\sqrt{n})^2$ with $K>1$ which is close to ours for $r\ll \min(q,n)$. For moderate to large $r$, we mention that our penalty~(\ref{condpen0}) can be significantly smaller than $r(\sqrt{q}+\sqrt{n})^2$, see Figure~\ref{fig1} below. 

\subsection*{Notations} 
All along the paper, we write $A^*$ for the adjoint of the matrix $A$ and $\sigma_{1}(A)\geq \sigma_{2}(A) \geq \ldots$ for its singular values ranked in a decreasing order. The Hilbert-Schmit norm of $A$ is denoted by $\|A\|=\mathrm{Tr}(A^*A)^{1/2}$ and the Ky-Fan $(2,r)$-norm by
$$\kf{A}{r}= \pa{\sum_{k=1}^r\sigma_{k}(A)^2}^{1/2}.$$
Finally, for a random variable $X$, we write $\E[X]^2$ for $\pa{\E[X]}^2$ to avoid multiple parentheses.

\section{A few facts on $\S_{q\times n}(r)$ and $X\widehat A_{r}$}
\subsection{Bounds on $\S_{q\times n}(r)$}
The expectation $\S_{q\times n}(r)=\E\cro{\kf{G_{q\times n}}{r}}$ can be evaluated numerically by Monte Carlo with a few lines of {\tt R}-code, see the Appendix. 
From a more theoretical point of view, we have the following bounds.
\begin{lemma} \label{thm singular}
 Assume that $q\leq n$.
Then for any $r\leq q$, we have $\S_{q\times n}(r)^2\geq r(n-1/q)$ and
$$\S_{q\times n}(r)^2
\leq \min\ac{r\,(\sqrt{n}+\sqrt{q})^2\,,nq-\sum_{k=r+1}^q(\sqrt{n}-\sqrt{k})^2,r+ \sum_{k=1}^r \pa{\sqrt{n}+\sqrt{q-k+1}}^2}.$$
When $q>n$ the same result holds with $q$ and $n$ switched. In particular, for $r=\min(n,q)$, we have
$$qn-1\leq \S_{q\times n}^2(\min(n,q))=\E\cro{\|G_{q\times n}\|}^2\leq qn.$$
\end{lemma}
The proof of the lemma is delayed to Section~\ref{sect-proof}. The map $r\to \S_{q\times n}(r)^2$ and the upper/lower bound of Lemma~1 are plotted in Figure~1 for $q=200$ and $n=200$ and 1000. We notice that the bound  $r\to \S_{q\times n}(r)^2\leq r(\sqrt{q}+\sqrt{n})^2$ looks sharp for small values of $r$, but it is quite loose for moderate to large values of $r$
\begin{figure}[h]
   \includegraphics[width=12cm,height=6.5cm]{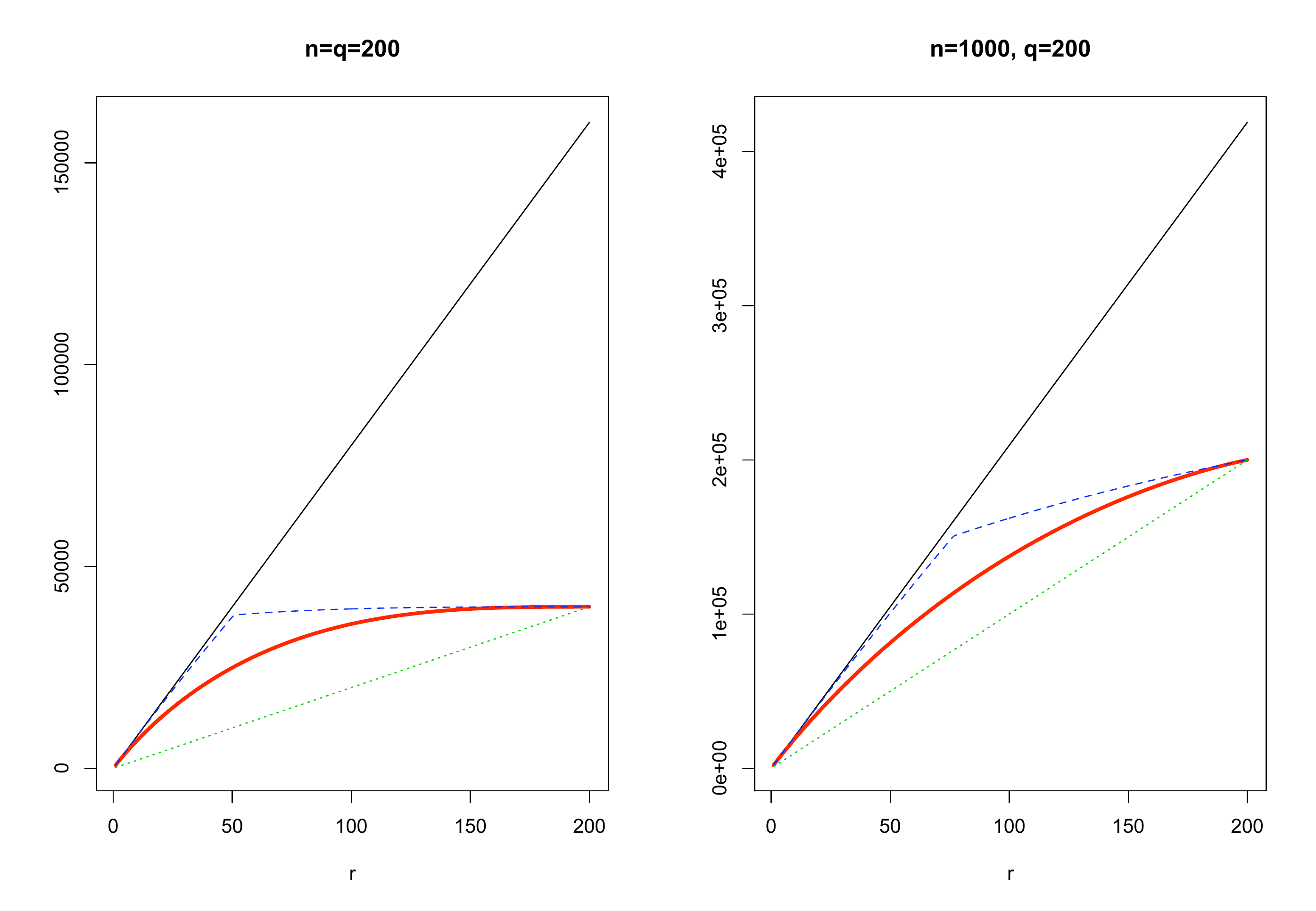}
      \caption{\label{fig1} In bold red $r \to \S_{q\times n}(r)^2$, in solid black $r\to r\,(\sqrt{n}+\sqrt{q})^2$, in dashed blue the upper-bound of Lemma~\ref{thm singular}, in dotted green the lower bound. Left: $q=n=200$. Right: $q=200$ and $n=1000$.}
\end{figure}

Finally, for large values of $q$ and $n$, asymptotics formulaes for $\S_{q\times n}(r)$ can be useful.
It is standard that when $n,q$ go to infinity with $q/n\to \beta\leq 1$, the empirical distribution of the eigenvalues of $n^{-1}G_{q\times n}G_{q\times n}^*$ converges almost surely to the  the Marchenko-Pastur distribution~\cite{MarchenkoPastur}, which has a density on $[(1-\sqrt{\beta})^2,(1+\sqrt{\beta})^2]$ given by
$$f_{\beta}(x)={1\over 2\pi\beta x}\, \sqrt{\big(x-(1-\sqrt{\beta})^2\big)\big((1+\sqrt{\beta})^2-x\big)}.$$
As a consequence, when $q$ and $n$ go to infinity with $q/n\to \beta\leq1$ and $r/q \to \alpha\leq1$, we have
 \begin{equation}\label{MP}
 \S_{q\times n}(r)^2\sim nq\int_{x_{\alpha}}^{(1+\sqrt{\beta})^2}xf_{\beta}(x)\, dx,
 \end{equation}
where
 $x_{\alpha}$ is defined by
$$\int_{x_{\alpha}}^{(1+\sqrt{\beta})^2}f_{\beta}(x)\, dx=\alpha.$$
Since the role of $q$ and $n$ is symmetric, the same result holds when $n/q\to \beta\leq1$ and $r/n \to \alpha\leq1$. 
This approximation~(\ref{MP}) can be evaluated efficiently (see the Appendix) and 
 it turns to be a very accurate approximation of $\S_{q\times n}(r)$ for $n,q$ large enough (say $nq>1000$).

\subsection{Computation of $X\hat A_{r}$}
Next lemma provides a useful formula for $X\hat A_{r}$.
\begin{lemma}\label{calcul}
Write $P$ for the projection matrix $P=X(X^*X)^+X^*$, with $(X^*X)^+$ the Moore-Penrose pseudo-inverse of $X^*X$. Then,  for any $r\leq q$ we have
$X\widehat A_{r}=(PY)_{r}$ where $(PY)_{r}$ minimizes $\|PY-B\|^2$ over the matrices $B$ of rank at most $r$.

As a consequence, writing $PY=U\Sigma V^*$ for  the singular value decomposition of $PY$, the matrix
$X\hat A_{r}$ is given by $X\hat A_{r}=U\Sigma_{r} V^*$, where $\Sigma_{r}$ is obtained from $\Sigma$ by setting $(\Sigma_{r})_{i,i}=0$ for $i\geq r+1$.
\end{lemma}
{\bf Proof of Lemma~\ref{calcul}.}
We note that  $\|PY-P(PY)_{r}\|^2\leq \|PY-(PY)_{r}\|^2$ and rank$(P(PY)_{r})\leq r$, so
$P(PY)_{r}=(PY)_{r}$. In particular, we have $(PY)_{r}=X\tilde A_{r}$, with $\tilde A_{r}=(X^*X)^+X^*(PY)_{r}$. Since the rank of $X\widehat A_{r}$ is also at most $r$, we have
\begin{eqnarray*}
\|Y-X\tilde A_{r}\|^2 &=& \|Y-PY\|^2+\|PY-(PY)_{r}\|^2\\
&\leq&  \|Y-PY\|^2+\|PY-X\widehat A_{r}\|^2=\|Y-X\widehat A_{r}\|^2.
\end{eqnarray*}
Since the rank of $\tilde A_{r}$ is not larger than $r$, we then have $\tilde A_{r}=\hat A_{r}$.

\section{The case of known variance}
In this section we revisit the results of Bunea {\it et al.}~\cite{BuneaSheWegkamp,Oct2010} for the case where $\sigma^2$ is known. 
This analysis will give us a benchmark  for the case of unknown variance.
Next theorem states an oracle inequality for the selection Criterion~(\ref{critvarianceconnue})
with penalty fulfilling $\pen_{\sigma^2}(r)\geq K \S_{q\times n}(r)^2$ for $K>1$. Later on, we will prove that the penalty $\pen_{\sigma^2}(r)= \S_{q\times n}(r)^2$ is minimal in some sense.
\begin{theorem}\label{thm oracle1}
Assume that for some $K>1$ we have 
\begin{equation}\label{condpen0}
\pen_{\sigma^2}(r)\geq K \S_{q\times n}(r)^2\quad \textrm{for all }\ \ r\leq \min(n,q). 
\end{equation}
Then, when $\hat r$ is selected by minimizing~(\ref{critvarianceconnue}) the estimator $\widehat A=\widehat A_{\hat r}$ satisfies
\begin{equation}\label{oracle1}
\E\cro{\|X\widehat A-XA\|^2} \leq c(K) \, \min_{r}\ac{\E\cro{\|XA-X\widehat A_{r}\|^2}+\pen_{\sigma^2}(r)\sigma^2+\sigma^2}
\end{equation}
for some  positive constant  $c(K)$ depending on $K$ only.
\end{theorem}
The risk bound~(\ref{oracle1}) ensures that the risk of the estimator $\widehat A$ is not larger (up to a constant) than the minimum over $r$ of the sum of the risk of the estimator $\widehat A_{r}$ plus the penalty term $\pen_{\sigma^2}(r)\sigma^2$. We will see below that this ensures that the estimator $\widehat A$ is adaptive minimax.\medskip

For $r\ll \min(n,q)$, the penalty $\pen_{\sigma^2}(r)= K \S_{q\times n}(r)^2$ is close to the penalty
$\pen'_{\sigma^2}=K(\sqrt{q}+\sqrt{n})^2r$ proposed by Bunea {\it et al.}~\cite{Oct2010}, but $\pen_{\sigma^2}(r)$ can be significantly smaller than $\pen'_{\sigma^2}(r)$ for moderate values of $r$, see Figure~1.
Next proposition shows that choosing a penalty  $\pen_{\sigma^2}(r)= K \S_{q \times n}(r)^2$ with $K<1$ can lead to a strong overfitting.
\begin{proposition}\label{propminimal1}
Assume that $A=0$ and that $\hat r$ is any minimizer of the Criterion~(\ref{critvarianceconnue}) with $\pen_{\sigma^2}(r)= K\S_{q \times n}(r)^2$ for some $K<1$. Then, setting $\alpha=1-\sqrt{(1+K)/2}>0$ we have
$$\p\pa{\hat r\geq {1-K\over 4}\times {nq-1\over (\sqrt{n}+\sqrt{q})^2}}\geq 1-e^{\alpha^2/2}\,{e^{-\alpha^2\max(n,q)/2} \over 1-e^{-\alpha^2\max(n,q)/2}}\,.$$
As a consequence,  the risk bound~(\ref{oracle1}) cannot hold when Condition~(\ref{condpen0}) is replaced by $\pen_{\sigma^2}(r)= K\S_{q \times n}(r)^2$ with $K<1$.
\end{proposition}
In the sense of Birg\'e and Massart~\cite{BirgeMassart}, the Condition~(\ref{condpen0}) is therefore minimal.

\subsection*{Minimax adaptation}
\begin{fact}
For any $\rho\in\,]0,1]$, there exists a constant $c_{\rho}>0$ such that for any integers $m,n,p$ larger than 2, any positive integer $q$ less than $\min(m,p)$ and any design matrix $X$ fulfilling
\begin{equation}\label{condminimax}
\sigma_{q}(X) \geq \rho\, \sigma_{1}(X),\quad\textrm{where}\ \ q=\textrm{rank}(X),
\end{equation}
we have
$$\inf_{\tilde A}\sup_{A\,:\, \mathrm{rank}(A)\leq r}\E\cro{\|X\tilde A-XA\|^2}\geq c_{\rho}(q+n)r\sigma^2,\quad\textrm{for all}\ \ r\leq \min(n,q).$$
\end{fact}
When $p\leq m$ and $q=p$, this minimax bound follows directly from Theorem~5 in Rohde and Tsybakov~\cite{RohdeTsybakov} as noticed by Bunea {\it et al.}, see~\cite{Oct2010} Section~2.3 Remark (ii) for a slightly different statment of this bound. We refer to Section~\ref{proof-minimax} for a proof of the general case (with possibly $q<p$ and/or $p>m$).

If we choose $\pen_{\sigma^2}(r)= K \S_{q\times n}(r)^2$ for some $K>1$, we have $\pen_{\sigma^2}(r)\leq 2Kr(q+n)$ according to Lemma~1. The risk bound~(\ref{oracle1}) then ensures that our estimator $\widehat A$ is adaptive minimax (as is the estimator proposed by Bunea {\it et al.}).

\section{The case of unknown variance}\label{sect-unknown}
We present now our main result which provides a selection criterion for the case where the variance $\sigma^2$ is unknown. For a given $r_{\max}\leq \min(n,q)$, we propose to select $\hat r\in\ac{1,\ldots,r_{\max}}$ by minimizing over $\ac{1,\ldots,r_{\max}}$  Criterion~(\ref{critunknown2}), namely
$$\crit(r)=\log(\|Y-X\widehat A_{r}\|^2)+\pen(r).$$
We note that the Criterion~(\ref{critunknown2}) is equivalent to the criterion
\begin{equation}\label{critunknown1}
\crit'(r)=\|Y-X\widehat A_{r}\|^2\pa{1+{\pen'(r)\over nm}},
\end{equation}
with $\pen'(r)=nm(e^{\pen(r)}-1)$.
This last criterion bears some similitude with the Criterion~(\ref{critvarianceconnue}). Indeed, the Criterion~(\ref{critunknown1}) can be written as
$$\|Y-X\widehat A_{r}\|^2+\pen'(r)\hat \sigma^2_{r},$$
with $\hat \sigma^2_{r}=\|Y-X\widehat A_{r}\|^2/(nm)$, which is the maximum likelihood estimator of $\sigma^2$ associated to $\widehat A_{r}$. To facilitate comparisons with the case of known variance, we will work henceforth with the Criterion~(\ref{critunknown1}).
Next theorem provides an upper bound for the risk of the estimator $X\widehat A_{\hat r}$.

\begin{theorem}\label{thm oracle2}
Assume that for some $K>1$ we have both
\begin{equation}\label{condrmax}
K\S_{q \times n}(r_{\max})^2+1<nm
\end{equation}
\begin{equation}\label{condpen1}
\textrm{and}\quad
\pen'(r)\geq {K\S_{q \times n}(r)^2\over 1-{1\over nm}(1+K\S_{q \times n}(r)^2)},\quad\textrm{for}\ r\leq r_{\max}.
\end{equation}
Then, when $\hat r$ is selected by minimizing~(\ref{critunknown1})  over $\ac{1,\ldots,r_{\max}}$, the estimator $\widehat A=\widehat A_{\hat r}$ satisfies
\begin{multline}\label{oracle2}
\E\cro{\|X\widehat A-XA\|^2}\\ \leq c(K)\, \min_{r\leq r_{\max}}\ac{\E\cro{\|X\widehat A_{r}-XA\|^2}\pa{1+{\pen'(r)\over nm}}+(\pen'(r)+1)\sigma^2}.
\end{multline}
for some constant $c(K)>0$ depending only on $K$.
\end{theorem}
Let us compare Theorem~\ref{thm oracle2} with Theorem~\ref{thm oracle1}. The two main differences lie in Condition~(\ref{condpen1}) and in the form of the risk bound~(\ref{oracle2}). Condition~(\ref{condpen1}) is more stringent than Condition~(\ref{condpen0}). More precisely, when $r$ is small compared to $q$ and $n$, both conditions are close, but when $r$ is of a size comparable to $q$ or $n$, Condition~(\ref{condpen1}) is much stronger than~(\ref{condpen0}). In the case where $m=q$, it even  enforces a blow up of the penalty $\pen'(r)$  when $r$ tends to $\min(n,m)$. This blow up is actually necessary to avoid overfitting since,  in this case,  the residual sum of squares $\|Y-X\widehat A_{r}\|^2$ tends to 0 when $r$ increases. The second major difference between Theorem~\ref{thm oracle2} and Theorem~\ref{thm oracle1} lies in the multiplicative factor $(1+\pen'(r)/nm)$ in the right-hand side of the risk bound~(\ref{oracle2}). 
Due to this term, the bound~(\ref{oracle2}) is not (strictly speaking) an oracle bound. To obtain an oracle bound, we have to add a condition on $r_{\max}$ to ensure that 
 $\pen'(r)\leq C nm$ for all $r\leq r_{\max}$. Next corollary provides such a condition.
 
 \begin{corollary}\label{cor} 
 Assume that
 \begin{equation}\label{rmax}
 r_{\max} \leq \alpha \, {nm-1\over K(\sqrt{q}+\sqrt{n})^2}\quad \textrm{for some }\ 0<\alpha<1,
 \end{equation}
and set
$$\pen(r) = -\log(1-K\S_{q \times n}(r)^2/ (nm-1))\quad \textrm{for some } K>1.$$

Then, there exists $c(K,\alpha)>0$ such that,  
 when $\hat r$ is selected by minimizing~(\ref{critunknown2})  over $\ac{1,\ldots,r_{\max}}$,
we have the oracle inequality
$$\E\cro{\|X\widehat A-XA\|^2}\\ \leq c(K,\alpha)\, \min_{r\leq r_{\max}}\ac{\E\cro{\|X\widehat A_{r}-XA\|^2}+r(\sqrt{n}+\sqrt{q})^2\sigma^2+\sigma^2}.$$
 \end{corollary}
In particular, the estimator $\widehat A$ is adaptive minimax up to the rank $r_{\max}$ specified by~(\ref{rmax}). In the worst case where $m=q$, Condition~(\ref{rmax}) requires that $r_{\max}$ remains smaller than a fraction of $\min(n,q)$. In the more favorable case where $m$ is larger than $4q$, Condition~(\ref{rmax}) can be met with $r_{\max}=\min(q,n)$ for suitable choices of $K$ and $\alpha$. 
 
 \begin{figure}[h]
   \includegraphics[width=12cm]{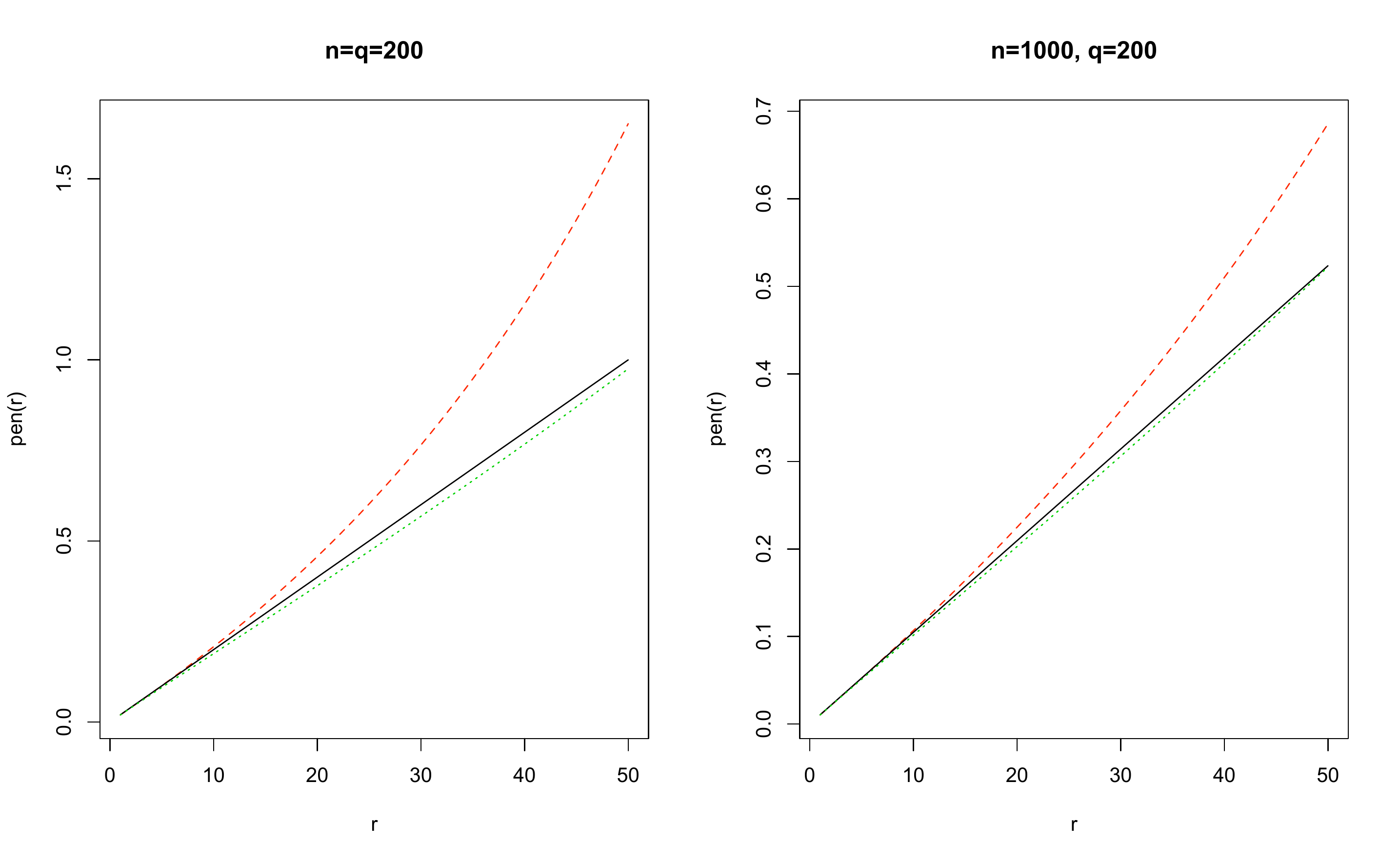}
      \caption{\label{fig2} In dotted green $\pen(r)=-\log(1-\S_{q \times n}(r)^2/ (nq-1))$, in solid black $\pen(r)=r\,(\sqrt{n}+\sqrt{q})^2/(nq)$, in dashed red $\pen'(r)/(nq)=\S_{q\times n}(r)^2/(nq-1-\S_{q\times n}(r)^2)$. Left : $q=n=200$. Right : $q=200$ and $n=1000$.}
\end{figure}

Let us discuss now in more details the Conditions~(\ref{condrmax}) and~(\ref{condpen1})  of Theorem~\ref{thm oracle2}.
We have 
$\S_{q\times n}(r)^2<{r(\sqrt{n}+\sqrt{q})^2}
$
so the Conditions~(\ref{condrmax}) and~(\ref{condpen1}) are satisfied as soon as 
$$r_{\max}\leq {nm-1\over K(\sqrt{n}+\sqrt{q})^2}
\quad\textrm{and}\quad
\ \pen'(r)\geq {Kr(\sqrt{n}+\sqrt{q})^2\over 1-{1\over nm}(1+ Kr(\sqrt{n}+\sqrt{q})^2)},\quad\textrm{for }r\leq r_{\max}.$$
In terms of the Criterion~(\ref{critunknown2}), the Condition~(\ref{condpen1}) reads 
$$\pen(r)\geq -\log(1-K\S_{q\times n}(r)^2/ (nm-1)).$$
When $\pen(r)$ is defined by taking equality in the above inequality, we have $\pen(r)\approx K r(\sqrt{n}+\sqrt{q})^2/(nm)$ for small values of $r$, see Figure~\ref{fig2}.

Finally, next proposition, shows that the Condition~(\ref{condpen1}) on $\pen'(r)$ is necessary to avoid overfitting.
\begin{proposition}\label{propminimal2}
Assume that $A=0$ and that $\hat r$ is any minimizer of Criterion~(\ref{critunknown1}) over $\ac{1,\ldots,\min(n,q)-1}$ with 
\begin{equation}\label{condminimal}
\pen'(r)={K\S_{q\times n}(r)^2\over 1-{K\over nm}\S_{q\times n}(r)^2}\quad\textrm{for some}\ K<1.
\end{equation}
Then, setting $\alpha=(1-K)/4>0$ we have
$$\p\pa{\hat r\geq {1-K\over 8}\times {nq-1\over (\sqrt{n}+\sqrt{q})^2}}\geq 1-2e^{\alpha^2/2}{e^{-\alpha^2\max(n,q)/2}\over 1-e^{-\alpha^2\max(n,q)/2}}\,.$$
\end{proposition}
As in Proposition~\ref{propminimal1},
a consequence of Proposition~\ref{propminimal2} is that Theorem~\ref{thm oracle2} cannot hold with Condition~(\ref{condpen1}) replaced by~(\ref{condminimal}). Condition~(\ref{condpen1}) is then minimal in the sense of Birg\'e and Massart~\cite{BirgeMassart}.

\section{Comments and extensions}\label{extensions}
\subsection{Link with PCA}
In the case where $X$ is the identity matrix, namely $Y=A+E$, Principal Component Analysis (PCA) is a popular technique to estimate $A$. The matrix $A$ is estimated by projecting the data $Y$ on the $r$ first  principal components, the number $r$ of components being chosen according to empirical or asymptotical criterions.

It turns out that the projection of the data  $Y$ on the $r$ first  principal components coincides with the estimator $\widehat A_{r}$. The criterion~(\ref{critunknown2}) then provides a theoretically grounded way to select the number $r$ of components.  Theorem~\ref{thm oracle2} ensures that the risk of  the final estimate $\widehat A_{\hat r}$ nearly 
achieves the minimum over $r$ of the risks  $\E\cro{\|\hat A_{r}-A\|^2}$.

\subsection{Sub-Gaussian errors}
We have considered for simplicity the case of Gaussian errors, but 
the results can be extended to the case where the entries $E_{i,j}$ are i.i.d sub-Gaussian. In this case, the matrix $PE$ will play the role of the matrix $G_{q\times n}$ in the Gaussian case. 
More precisely, combining recent results of Rudelson and Vershynin~\cite{RudelsonVershynin} and Bunea {\it et al.}~\cite{BuneaSheWegkamp}  on sub-Gaussian random matrices, with concentration inequality for sub-Gaussian random variables~\cite{ledoux} enables to prove an analog of 
Lemma~\ref{thm singular} for $\E\cro{\kf{PE}{r}}^2$ (with different constants). Then, the proof of  Theorem~\ref{thm oracle1} and Theorem~\ref{thm oracle2} can be easily adapted, replacing the Condition~(\ref{condpen0}) by
$$\pen(r)\geq K \E\cro{\kf{PE}{r}}^2, \quad\textrm{for }r\leq\min(q,n),$$
and the Conditions~(\ref{condrmax}) and~(\ref{condpen1}) by 
$K\E\cro{\kf{PE}{r_{\max}}}^2<\E[\|E\|]^2$ and
$$\pen'(r)\geq {K\E\cro{\kf{PE}{r}}^2\over 1-K\E\cro{\kf{PE}{r}}^2/\,\E\cro{\|E\|}^2},\quad\textrm{for}\ r\leq r_{\max}.$$
Analogs of Proposition 1 and 2 also hold with different constants.

\subsection{Selecting among arbitrary estimators}
Our theory provides a procedure to select among the family of estimators $\{\widehat A_{r},\ r\leq r_{\max}\}$. It turns out that it can be extended to arbitrary (finite) families of estimators $\{A_{\lambda},\ \lambda\in\Lambda\}$ such as the nuclear norm penalized estimator family $\{\widehat A_{\lambda}^{\ell_{1}}, \lambda\in\Lambda\}$. The most straightforward way is to replace everywhere $\widehat A_{r}$ by $\widehat A_{\lambda}$ and $\pen(r)$ by $\underline{\pen}(\lambda)$, with
$\underline{\pen}(\lambda)=\pen(\textrm{rank}(\widehat A_{\lambda})).$
In the spirit of Baraud {\it et al.}~\cite{linselect}, we may also consider more refined criterions such as 
$$\crit_{\alpha}(\lambda)=\min_{r\leq r_{\max}} \ac{(\|Y-X\widehat A_{\lambda,r}\|^2+\|X\widehat A_{\lambda}-X\widehat A_{\lambda,r}\|^2)\pa{1+{\pen'(r)\over nm}}},$$
where $\alpha>0$ and $\widehat A_{\lambda,r}$ minimizes $\|B-\widehat A_{\lambda}\|$ over the matrices $B$ of rank at most $r$. Analogs of Theorem~\ref{thm oracle2} can be derived for such criterions, but we will not pursue in that direction.

\section{Numerical experiments}
We perform numerical experiments on synthetic data in  two different settings. In the first experiment, we consider a favorable setting where the sample size $m$ is large compared to the number $p$ of covariables. In the second experiment, we consider a more challenging setting where the sample size $m$ is small compared to $p$. The objectives of our experiments are mainly:
\begin{itemize}
\item to give an example of implementation of our procedure,
\item to demonstrate that it can handle high-dimensional settings.
\end{itemize}
  
\subsection*{Simulation setting}
Our experiments are inspired by those of  Bunea {\it et al.}~\cite{BuneaSheWegkamp, Oct2010}, the main difference is that we work in higher dimensions. 
The simulation design is the following. The rows of the matrix $X$   are drawn independently according to a centered Gaussian distribution with covariance matrix  $\Sigma_{i,j}=\rho^{|i-j|}$, $\rho>0$. For a positive $b$, the matrix $A$ is given by $A=bB_{p\times r}B_{r \times n}$, where the entries of the $B$ matrices are i.i.d. standard Gaussian. 
For $r\leq \min(n,p)$, the rank of the matrix $A$ is then $r$ with probability one and the rank of $X$ is $\min(m,p)$ a.s.

\subsubsection*{Experiment 1:} in the first experiment, we consider a case where the sample size 
$m=400$ is large compared to the number $p=100$ of covariables and $n=100$. The other parameters are $r=40$,   $\rho$ varies in $\{0.1, 0.5, 0.9\}$ and $b$ varies in $\{0.025, 0.05, 0.075, 0.1\}$.
This experiment is actually the same as the Experiment~1 in~\cite{Oct2010}, except that we have multiplied $m,\ p,\ n,\ r$  by four and adjusted the values of $b$.

\subsubsection*{Experiment 2:} the second experiment is much more challenging since the sample size $m=100$ is small compared to 
 the number $p=500$ of covariables and $n=500$. Furthermore, the rank $q$ of $X$ equals $m$, which is the least-favorable case for estimating the variance. For the other parameters, we set $r=20$, $\rho$ varies in $\{0.1, 0.5, 0.9\}$ and $b$ varies in $\{0.005, 0.01, 0.015, 0.02\}$. 

\subsection*{Estimators}
For $K>1$, we write $\KF[K]$ for the estimator $\widehat A_{\hat r}$ with $\hat r$ selected by  the  Criterion~(\ref{critunknown1}) with
$$\pen'(r) = {K\S_{q\times n}(r)^2\over 1-{1\over nm}(1+K\S_{q\times n}(r)^2)}$$
(the notation $\KF$ refers to the Ky-Fan norms involved in $\S_{q\times n}(r)$).
\smallskip

For $\lambda>0$, we write $\RSC[\lambda]$ for the estimator $\widehat A_{\hat r}$ with $\hat r$ selected by  the criterion introduced by Bunea, She and Wegkamp~\cite{BuneaSheWegkamp}
$$\crit_{\lambda}(r)=\|Y-X\widehat A_{r}\|^2+\lambda (n+\textrm{rank}(X))  r.$$
Bunea {\it et al.}~\cite{Oct2010} proposes to use $\lambda=K\hat \sigma^2$ with $K\geq 2$ and
$$\hat \sigma^2={\|Y-PY\|^2\over mn-n\textrm{rank}(X)},\quad\textrm{with $P$ the projector onto the range of }X.$$
We denote by $\RSCI[K]$ the resulting estimator $\RSC[K\hat\sigma^2]$. 
\smallskip

Both procedures $\KF$ and $\RSCI$ depend on a tuning parameter $K$.  There  is no reason for the existence of a universal "optimal" constant $K$.
Nevertheless, Birg\'e and Massart~\cite{BirgeMassart} suggest to penalize by twice the minimal penalty, which corresponds to the choice $K=2$ for $\KF$. The value $K=2$ is also the value recommended by Bunea {\it et al.}~\cite{Oct2010}  Section~4 for the $\RSCI$ (see the "adaptive tuning parameter" $\mu_{adapt}$). Another classical approach for choosing a tuning parameter is Cross-Validation : for example, $K$ can be selected among a small grid of values between 1 and 3 by $V$-Fold CV. We emphasize that there is no theoretical justification that Cross-Validation will choose the {\it best} value for $K$. Nevertheless, for {\it each} value $K$ in the grid, the estimators $\KF[K]$ and $\RSCI[K]$ fulfills an oracle inequality with large probability, so the estimators with $K$ chosen by CV will also fulfills an oracle inequality with large probability (as long as the size of the grid remains small). We will write $\KF[K\!=\!CV]$ and $\RSCI[K\!=\!CV]$ for the estimators $\KF$ and $\RSCI$ with $K$ selected by $10$-fold Cross-Validation.

Finally, in Experiment~2 the estimator $\RSCI$ cannot be implemented since rank$(X)=m$. Yet, it is still possible to implement the procedure $\RSC[\lambda]$ and select $\lambda>0$ among a large grid of values  by 10-fold Cross-Validation, even if in this case there is no theoretical control on the risk of the resulting estimator $\RSC[\lambda\!=\!CV]$. We will use this estimator as a benchmark in Experiment~2.

\subsection*{Results}
\begin{figure}[h!]
   \includegraphics[width=\textwidth, height=0.9\textheight]{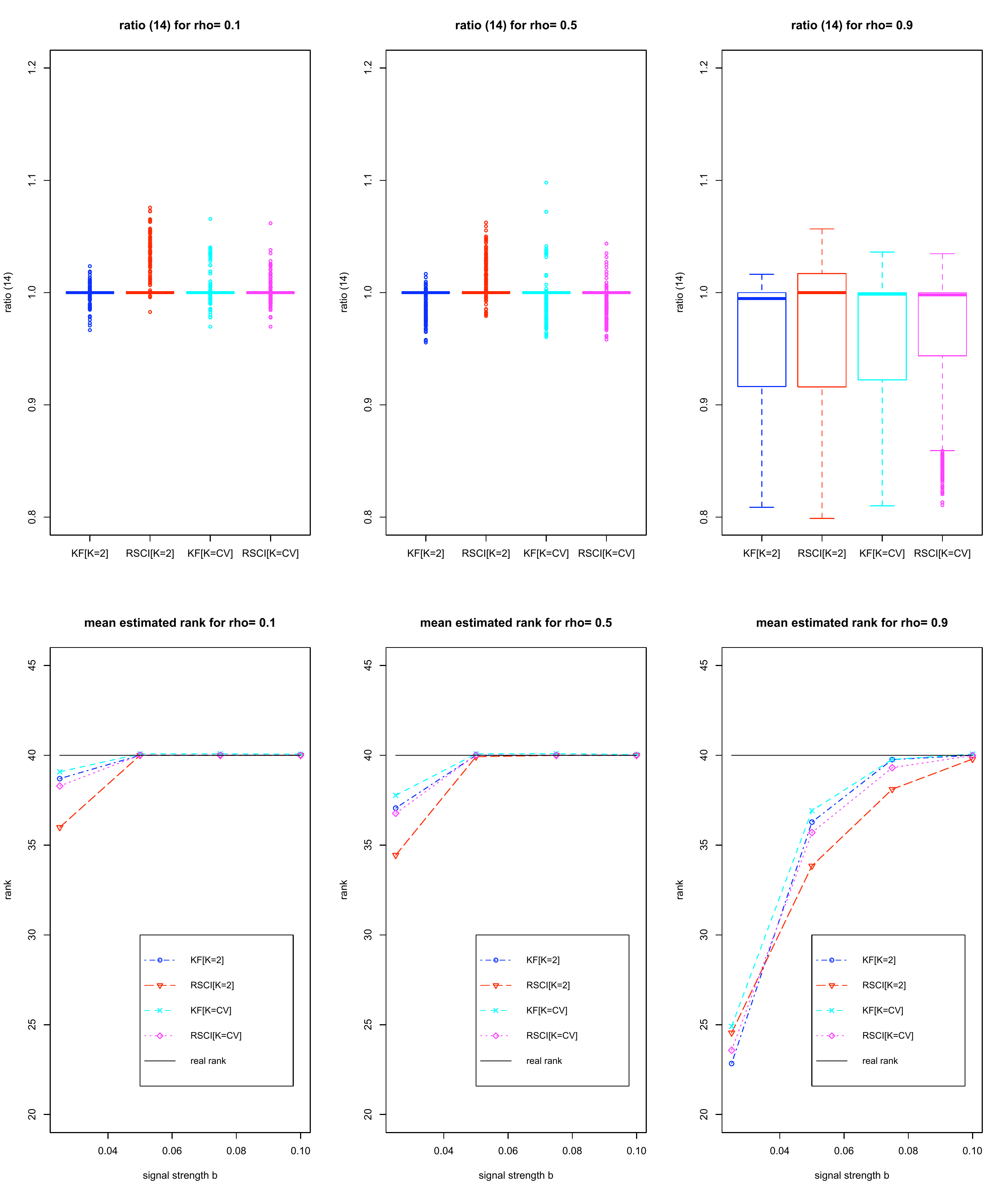} 
      \caption{\label{plot1}{\bf Experiment 1.} Left to right : $\rho=0.1,\,0.5,\,0.9$. Top : boxplots of the ratio~(\ref{boxplots}) for $\KF[K\!=\!2]$,  $\RSCI[K\!=\!2]$, $\KF[K\!=\!CV]$ and $\RSCI[K\!=\!CV]$. Bottom : mean estimated rank $\E[\hat r]$  for each estimator and each value of $b$.}
\end{figure}

The results of the first experiment are reported in Figure~\ref{plot1} and those of the second experiment in Figure~\ref{plot2}. 
The boxplots of the first line compare the performances of estimators $\KF$ and $\RSCI$ to that of the estimator $X\widehat A_{r}$ that we would use if we knew that the rank of A is $r$. The boxplots give for each value of $\rho$ the distribution of the ratios
\begin{equation}\label{boxplots}
{\|XA-X\widehat A\|^2\over\|XA-X\widehat A_{r}\|^2}\  \wedge\ 10
\end{equation}
for $\widehat A$ given by $\KF[K\!=\!2]$, $\RSCI[K\!=\!2]$, $\KF[K\!=\!CV]$ and $\RSCI[K\!=\!CV]$ in the first experiment and
by $\KF[K\!=\!2]$, $\KF[K\!=\!CV]$, and $\RSC[\lambda\!=\!CV]$ in the second experiment.
The ratios~(\ref{boxplots}) are truncated to 10 for a better visualization.
Finally, we plot in the second line the mean  estimated rank $\E[\hat r]$  for each estimator and each value of $b$ and $\rho$.

\subsubsection*{Experiment~1 (large sample size)} All estimators $\KF[K\!=\!2]$, $\RSCI[K\!=\!2]$, $\KF[K\!=\!CV]$ and $\RSCI[K\!=\!CV]$ perform very similarly and  almost all the ratios~(\ref{boxplots}) are equal to~1.

\subsubsection*{Experiment~2 (small sample size)} 
\begin{figure}[h!]
   \includegraphics[width=\textwidth,height=0.9\textheight]{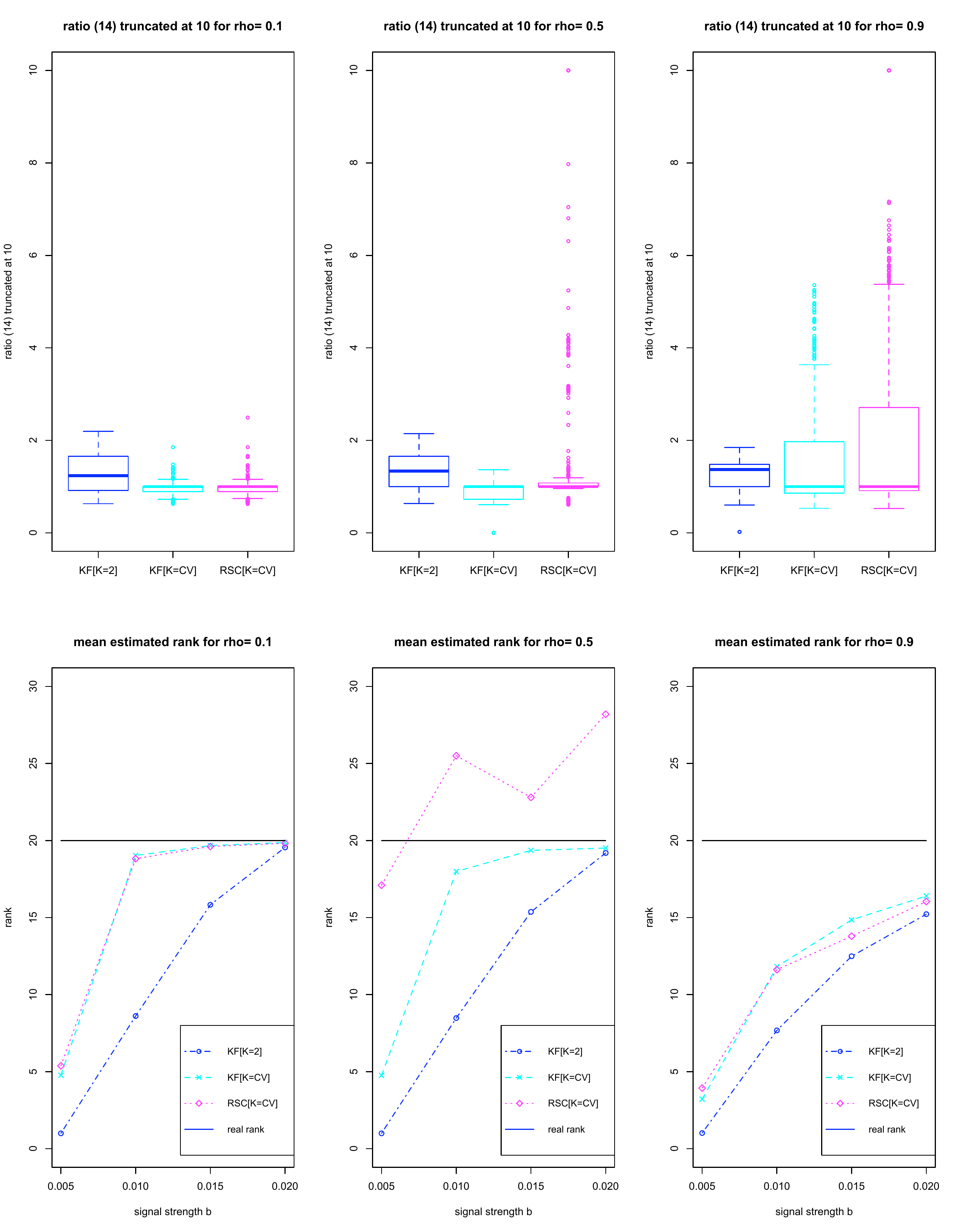} 
      \caption{\label{plot2}{\bf Experiment 2.}  Left to right : $\rho=0.1,\,0.5,\,0.9$. Top : boxplots of the ratio~(\ref{boxplots}) truncated at 10 for $\KF[K\!=\!2]$,  $\KF[K\!=\!CV]$ and $\RSC[\lambda\!=\!CV]$. Bottom : mean estimated rank $\E[\hat r]$  for each estimator and each value of $b$.}
\end{figure}

The estimator $\KF[K\!=\!CV]$ has global good performances, with a median ratio~(\ref{boxplots}) around~1, but the ratio~(\ref{boxplots}) can be as high as 5 in some examples for $\rho=0.9$.  In contrast,  the estimator $\KF[K\!=\!2]$ is very stable but it has a median value significantly above the other methods.
Finally, the performances of the estimator $\RSC[\lambda\!=\!CV]$ are very contrasted. For small correlation ($\rho=0.1$), its performances are similar to that of $\KF[K\!=\!CV]$. For $\rho=0.5$ or $\rho=0.9$, it has very good performances most of the time  (similar to $\KF[K\!=\!CV]$) but  it completely fails on a small fraction of examples. For example, for $\rho=0.9$, it has a ratio~(\ref{boxplots}) smaller than $7$ in 80\% of the examples (with a median value close to 1), but in
20\% of the examples, it completely fails and the ratio $\|XA-X\widehat A\|^2/\|XA-X\widehat A_{r}\|^2$ for $\RSC[\lambda\!=\!CV]$ can be has high as $10^{13}$ (these values do not appear in Figure~\ref{plot2} since we have truncated the ratio~(\ref{boxplots}) to 10 to avoid a complete shrinkage of the boxplots). We recall, that there exists no risk bound for the estimator $\RSC[\lambda\!=\!CV]$, so these results are not in contradiction with any theory.

Finally,  we emphasize that no conclusion should be drawn from these two single experiments about the superiority of one procedure over the others.


\section{Proof of the mains results}\label{sect-proof}
\subsection{Proof of Lemma~\ref{thm singular}}
 For notational simplicity we write $G=G_{q\times n}$.
The case $r=1$ follows from  Slepian's Lemma, see Davidson and Szarek~\cite{DavidsonSzarek} Chapter~8. For $r>1$, we note that
$$\E\cro{\kf{G}{r}}^2\leq \min\ac{r\,\E\cro{\kf{G}{1}}^2\,,\sum_{k=1}^r \E[\sigma_{k}^2(G)]}.$$
The first upper bound $\S_{q\times n}(r)^2\leq r(\sqrt{n}+\sqrt{q})^2$ follows. For the second upper bound, 
we note that
$$\sum_{k=1}^r \E[\sigma_{k}^2(G)]\leq \E[\|G\|^2]-\sum_{k=r+1}^q \E[\sigma_{k}(G)]^2.$$
The  interlacing inequalities~\cite{HornJohnson} ensure that  $\sigma_{k}(G_{k}')\leq \sigma_{k}(G)$ 
where $G_{k}'$ is the matrix made of the $k$ first rows of $G$.
The second bound then follows from $\E[\sigma_{k}(G_{k}')]\geq \sqrt{n}-\sqrt{k}$, see~\cite{DavidsonSzarek}.

Let us turn to the third bound. The map $G \to \sigma_{k}(G)$ is 1-Lipschitz so, writing $M_{k}$ for the median of $\sigma_{k}$, the concentration inequality for Gaussian random variables ensures that $(M_{k}-\sigma_{k}(G))_{+}\leq \xi_{+}$ and $ (\sigma_{k}(G)-M_{k})_{+}\leq\xi_{+}'$ where $\xi_{+}$ and  $\xi_{+}'$ are the positive part of two standard Gaussian random variables. 
As a consequence we have
\begin{eqnarray*}
\E[\sigma_{k}^2(G)]-\E[\sigma_{k}(G)]^2+(M_{k}-\E[\sigma_{k}(G)])^2&=&\E\cro{(\sigma_{k}(G)-M_{k})_{+}^2}+\E\cro{(M_{k}-\sigma_{k}(G))_{+}^2}\nonumber \\
&\leq& \E[\xi_{+}^{'2}]+\E[\xi_{+}^{2}]=1,\nonumber
\end{eqnarray*}
and thus $\E[\sigma_{k}^2(G)]\leq\E[\sigma_{k}(G)]^2+1$.

Furthermore, the interlacing inequalities~\cite{HornJohnson} ensure that $\sigma_{k}(G)\leq \sigma_{1}(G_{q-k+1}')$. We can then bound $\E[\sigma_{k}(G)]$ by
$$\E[\sigma_{k}(G)]\leq \sqrt{n}+\sqrt{q-k+1}$$
which leads to the last upper bound.

For the lower bound, we start from $\kf{G}{r}^2\geq \|G\|^2r/q$ (sum of a decreasing sequence) and use again the Gaussian concentration inequality to get
$$\E[\|G\|^2]-1=nq-1\leq \E[\|G\|]^2$$
and concludes that $r(nq-1)/q\leq \E\cro{\kf{G}{r}}^2=\S_{q\times n}(r)^2$.

\subsection{A technical lemma}
Next lemma provides a control of the size of the scalar product $<E,X\widehat A_{k}-XA_{r}>$ which will be useful for the proofs of Theorem~\ref{thm oracle1} and Theorem~\ref{thm oracle2}.
\begin{lemma}\label{tech}
Fix $r\leq \min(n,q)$ and write  $A_{r}$ for the best approximation of $A$ with rank at most $r$. Then, there exists  a random variable $U_{r}$ such that  $\E(U_{r})\leq r\min(n,q)$ and for any $\eta>0$ and all $k\leq \min(n,q)$
\begin{eqnarray} \label{ineq1}
2 \sigma\, |<E,X\widehat A_{k}-XA_{r}>|
&\leq & {1\over 1+\eta}\|X\widehat A_{k}-XA\|^2+{1+1/\eta\over (1+\eta)^2}\|XA-XA_{r}\|^2\\ \nonumber
& & \ +\, (1+\eta)^2(1+1/\eta)\sigma^2U_{r}+ (1+\eta)^3\sigma^2\kf{PE}{k}^2
\end{eqnarray}
where $P=X(X^*X)^+X^*$ is as in Lemma~\ref{thm singular}.
\end{lemma}
Iterating twice the inequality $2ab\leq a^2/c+cb^2$ gives
\begin{multline*}
2 \sigma\, |<E,X\widehat A_{k}-XA_{r}>|
\\ \leq {1\over 1+\eta}\|X\widehat A_{k}-XA\|^2+{1+1/\eta\over (1+\eta)^2}\|XA-XA_{r}\|^2+(1+\eta)^2\sigma^2{<E,X\widehat A_{k}-XA_{r}>^2\over \|X\widehat A_{k}-XA_{r}\|^2}.
\end{multline*}
We write $XA_{r}=U\Gamma_{r}V^*$ for the singular value decomposition of $XA_{r}$, with the convention that the diagonal entries of $\Gamma_{r}$ are decreasing. Since the rank of $XA_{r}$ is upper bounded by the rank of $A_{r}$, the $m\times n$ diagonal matrix $\Gamma_{r}$ has at most $r$ non zeros elements. 
Assume first that $n\leq q$.
Denoting by  $I_{r}$  the $m\times m$ diagonal matrix with $(I_{r})_{i,i}=1$ if $i\leq r$ and $(I_{r})_{i,i}=0$ if $i>r$ and writing $I_{-r}=I-I_{r}$ and $\widehat B_{k}=U^*X\widehat A_{k}V$, we have
\begin{eqnarray*}
{<E,X\widehat A_{k}-XA_{r}>^2\over \|X\widehat A_{k}-XA_{r}\|^2}&=&{<U^*PEV,\widehat B_{k}-\Gamma_{r}>^2\over \|\widehat B_{k}-\Gamma_{r}\|^2}\\
&=&{\pa{<U^*PEV,I_{r}(\widehat B_{k}-\Gamma_{r})>+<U^*PEV,I_{-r}\widehat B_{k}>}^2\over \| I_{r}(\widehat B_{k}-\Gamma_{r})\|^2+\|I_{-r}\widehat B_{k}\|^2}\\
&\leq& (1+\eta^{-1}) {<U^*PEV,I_{r}(\widehat B_{k}-\Gamma_{r})>^2\over \| I_{r}(\widehat B_{k}-\Gamma_{r})\|^2} 
+(1+\eta){<U^*PEV,I_{-r}\widehat B_{k}>^2\over \|I_{-r}\widehat B_{k}\|^2}.
\end{eqnarray*}
The first term is upper bounded by
 $${<U^*PEV,I_{r}(\widehat B_{k}-\Gamma_{r})>^2\over \| I_{r}(\widehat B_{k}-\Gamma_{r})\|^2} \leq \| I_{r}U^*PEV\|^2=U_{r}$$
 and the expected value of the right-hand side fulfills
 $$\E(U_{r})=n\|I_{r}U^*PU\|^2=n\|U^*PUI_{r}\|^2\leq nr.$$
Since the rank of $I_{-r}\widehat B_{k}$ is at most $k$, the second term can be bounded by
$${<U^*PEV,I_{-r}\widehat B_{k}>^2\over \|I_{-r}\widehat B_{k}\|^2}\leq \sup_{\mathrm{rank}(B)\leq k}{<U^*PEV,B>^2\over \|B\|^2}=\kf{U^*PEV}{k}^2=\kf{PE}{k}^2.$$
Putting pieces together gives~(\ref{ineq1}) for $n\leq q$. The case  $n>q$ can be treated in the same way, starting from 
$$\widehat B_{k}-\Gamma_{r}=(\widehat B_{k}-\Gamma_{r})I_{r}+\widehat B_{k}I_{-r}$$
with $I_{r}$ and $I_{-r}$  two $n\times n$ diagonal matrices defined as above.

\subsection{Proof of Theorem~\ref{thm oracle1}}\label{proof-thm1}
The inequality $\crit_{\sigma^2}(\hat r)\leq \crit_{\sigma^2}(r)$ gives 
\begin{equation}\label{ineq0}
\|X\widehat A-XA\|^2\leq \|X\widehat A_{r}-XA\|^2+\pen_{\sigma^2}(r)\sigma^2+2\sigma<E,X\widehat A-X\widehat A_{r}>-\pen_{\sigma^2}(\hat r) \sigma^2.
\end{equation}
Combining this inequality with Inequality~(\ref{ineq1}) of Lemma~\ref{tech} with $\eta=((1+K)/2)^{1/3}-1>0$, we obtain
 \begin{eqnarray*}
 {\eta\over 1+\eta} \|X\widehat A-XA\|^2&\leq& \|X\widehat A_{r}-XA\|^2+2{1+1/\eta \over (1+\eta)^2}\|XA-XA_{r}\|^2+2\pen_{\sigma^2}(r)\sigma^2\\
 &&\ +2(1+\eta)^2(1+\eta^{-1})\sigma^2U_{r}+{K+1\over 2}\sigma^2\kf{PE}{r}^2-\pen_{\sigma^2}(r)\sigma^2\\
&& \ +\sigma^2\sum_{k=1}^{\min(n,q)}\pa{{K+1\over 2}\,\kf{PE}{k}^2-\pen_{\sigma^2}(k)}_{+}.
\end{eqnarray*}
The map $E \to \kf{PE}{k}$ is 1-Lipschitz and convex, so there exists a standard Gaussian random variable $\xi$ such that 
$\kf{PE}{k}\leq\E[\kf{PE}{k}]+\xi_{+}$ and then
\begin{eqnarray*}
\E\pa{{K+1\over 2}\,\kf{PE}{k}^2-\pen(k)}_{+}&\leq& {1+K\over 2}\E\pa{\xi_{+}^2+2\xi_{+}\E[\kf{PE}{k}]-{K-1\over K+1}\E[\kf{PE}{k}]^2}_{+}\\
&\leq& c_{1}(K) \exp(-c_{2}(K) \E[\kf{PE}{k}]^2).
\end{eqnarray*}
Since $\kf{PE}{k}$ is distributed as $\kf{G_{q\times n}}{k}$, 
 Lemma~\ref{thm singular} gives that  $\E[\kf{PE}{k}]^2\geq k \max(n,q)-1$
and the series 
$$\sum_{k=1}^{\min(n,q)}\E\pa{{K+1\over 2}\,\kf{PE}{k}^2-\pen(k)}_{+}$$
can be upper-bounded by 
$c_{1}(K)e^{c_{2}(K)}\pa{1-e^{-c_{2}(K)}}^{-1}e^{-c_{2}(K)\max(n,q)}$.
Finally, $\E[U_{r}]\leq r\min(n,q)$ is bounded by  $1+\pen(r)$ and $\|XA-XA_{r}\|^2$ is smaller than $\E\cro{\|XA-X\widehat A_{r}\|^2}$, so there exists some constant $c(K)>0$ such that~(\ref{oracle1}) holds.
%

\subsection{Proof of  Theorem~\ref{thm oracle2}.}
To simplify the formulaes, we will note $\penn(r)=\pen'(r)/(nm)$.
The inequality $\crit'(\hat r)\leq \crit'(r)$ gives
\begin{eqnarray*}
\lefteqn{\|X\widehat A-XA\|^2(1+\penn(\hat r))}\\ &\leq& \|Y-X\widehat A_{r}\|^2-\sigma^2(1+\penn(r))\|E\|^2+\penn(r)\|Y-X\widehat A_{r}\|^2+\penn(r)\|E\|^2\sigma^2\\
&&+2(1+\penn(\hat r))\sigma<E,X\widehat A-XA>-\penn(\hat r)\|E\|^2\sigma^2\\
&\leq& \pa{2\sigma<E,XA_{r}-X\widehat A_{r}>-\penn(r)\sigma^2\|E\|^2}_{+}+(1+2\penn(r))\|XA-X\widehat A_{r}\|^2+3\penn(r)\|E\|^2\\
&&+(1+\penn(\hat r))\pa{2\sigma<E,X\widehat A-XA_{r}>-{\penn(\hat r)\over 1+\penn(\hat r)}\,\|E\|^2\sigma^2}_{+}+2\sigma\penn(\hat r)<E,XA_{r}-XA>.
\end{eqnarray*}
Dividing both side by $1+\penn(\hat r)$, we obtain
$$\|X\widehat A-XA\|^2\leq (1+2\penn(r))\|XA-X\widehat A_{r}\|^2+3\penn(r)\|E\|^2+2\sigma|<E,XA-XA_{r}>|+\Delta_{r}+\Delta_{\hat r}$$
where
$$\Delta_{k}=\pa{2\sigma|<E,X\widehat A_{k}-XA_{r}>|-{\penn(k)\over 1+\penn(k)}\,\|E\|^2\sigma^2}_{+}.$$
We first  note that $\E[\|E\|^2]=nm$ and $2\sigma\E[|<E,XA-XA_{r}>|]\leq \sigma^2+\|XA-XA_{r}\|^2$. 
Then, combining Lemma~\ref{tech} with $\eta=(K^{1/6}-1)$ and the following lemma with $\delta=\eta$ gives
$$\E\cro{\|X\widehat A-XA\|^2}\leq c(K) \pa{\E\cro{\|XA-X\widehat A_{r}\|^2}(1+\penn(r))+(1+nm\penn(r))\sigma^2},$$
for some $c(K)>0$.

\begin{lemma}
Write $P$ for the projection matrix $P=X(X^*X)^+X^*$, with $(X^*X)^+$ the Moore-Penrose pseudo-inverse of $X^*X$.
For any $\delta>0$ and $r\leq \min(n,q)$ such that $(1+\delta)\E\cro{\kf{PE}{r}}\leq \sqrt{nm-1}$, we have
\begin{equation}\label{controle}
\E\cro{(\kf{PE}{r}^2-(1+\delta)^3\E\cro{\kf{PE}{r}}^2\|E\|^2/(nm-1))_{+}}\leq 4\pa{1+1/\delta}e^{\delta^2/4}e^{-\delta^2r\max(n,q)/4}.
\end{equation}
As a consequence, we have
\begin{multline*}\E\cro{\sup_{r\leq r_{\max}}\pa{\kf{PE}{r}^2-(1+\delta)^3\S_{q\times n}(r)^2\|E\|^2/(nm-1)}_{+}}\\
\leq 4\pa{1+1/\delta}e^{\delta^2/4}{e^{-\delta^2\max(n,q)/4}\over1-e^{-\delta^2\max(n,q)/4}}\,.
\end{multline*}
\end{lemma}
{\bf Proof of the Lemma.}

Writing $t=(1+\delta)\E[\kf{PE}{r}]/\E[\|E\|]\leq 1$,  the map $E\to \kf{PE}{r}-t\|E\|$ is $\sqrt{2}$-Lipschitz. Gaussian concentration inequality then ensures that 
\begin{eqnarray*}
\kf{PE}{r}&\leq& t\|E\|+\E[\kf{PE}{r}-t\|E\|]+2\sqrt{\xi}\\
&\leq& t\|E\|+\pa{2\sqrt{\xi}-\delta\,\E[\kf{PE}{r}]}_{+},
\end{eqnarray*}
with $\xi$  a standard exponential random variable. We then get that
$$\kf{PE}{r}^2\leq (1+\delta)t^2\|E\|^2+4(1+1/\delta)\pa{\sqrt{\xi}-\delta\,\E[\kf{PE}{r}]/2}_{+}^2$$
and 
\begin{eqnarray*}
\E\cro{(\kf{PE}{r}^2-(1+\delta)t^2\|E\|^2)_{+}}&\leq&4(1+1/\delta)\E\cro{\pa{\sqrt{\xi}-\delta\,\E[\kf{PE}{r}]/2}_{+}^2}\\
&\leq&  4\pa{1+1/\delta}e^{-\delta^2\E[\kf{PE}{r}]^2/4}.\\
\end{eqnarray*}
The bound~(\ref{controle}) then follows from $\E[\kf{PE}{r}]^2\geq r\max(n,q)-1$ and $\E[\|E\|]^2\geq nm-1$.\hfill$\square$


\subsection{Proof Proposition~\ref{propminimal1}} For simplicity we consider first the case where $m=q$. With no loss of generality, we can also assume that $\sigma^2=1$.
We set
$$\Omega_{0}=\ac{\|E\|\geq (1-\alpha) \E[\|E\|]}\bigcap_{r=1}^{\min(n,m)} \ac{\kf{E}{r}\leq (1+\alpha)\E[\kf{E}{r}]}.$$
According to the Gaussian concentration inequality we have
\begin{eqnarray*}
\p(\Omega_{0})&\geq& 1-\sum_{r=1}^{\min(n,m)}e^{-\alpha^2\E[\kf{E}{r}]^2/2}\\
&\geq& 1-e^{\alpha^2/2}\sum_{r=1}^{\min(n,m)}e^{-\alpha^2r\max(n,m)/2}
\end{eqnarray*}
where the last bound follows from Lemma~\ref{thm singular}.
Furthermore, Lemma~\ref{calcul} gives that $X\hat A_{r}=Y_{r}(=E_{r})$, where $Y_{r}$ is the matrix $M$ minimizing $\|Y-M\|^2$ over the matrices of rank at most $r$.
As a consequence, 
writing $m^*=\min(n,m)$, we have on $\Omega_{0}$
\begin{eqnarray*}
\crit_{\sigma^2}(m^*)-\crit_{\sigma^2}(r)&=&K\E[\kf{E}{m^*}]^2-(\kf{E}{m^*}^2-\kf{E}{r}^2)-K\E[\kf{E}{r}]^2\\
&\leq& \pa{(1+\alpha)^2-K}\E[\kf{E}{r}]^2-\pa{(1-\alpha)^2-K}\E[\kf{E}{m^*}]^2\\
&\leq& 2\E[\kf{E}{r}]^2-{1-K\over 2}\E[\kf{E}{m^*}]^2\\
&<& 2r(\sqrt{n}+\sqrt{m})^2-{1-K\over 2}(nm-1).
\end{eqnarray*}
We then conclude that on $\Omega_{0}$ we have $\hat r\geq {1-K\over 4}\times {nm-1\over (\sqrt{n}+\sqrt{m})^2}$.

Let $r^*$ be the smaller integer larger than ${1-K\over 4}\times {nm-1\over (\sqrt{n}+\sqrt{m})^2}$.
Since $\|X\widehat A-XA\|^2=\kf{E}{\hat r}^2$, we have
\begin{eqnarray*}
\E\cro{\|X\widehat A-XA\|^2}&\geq& \E\cro{\kf{E}{r^*}^2\mathbf{1}_{\hat r\geq r^*}}\\
&\geq& (1-\alpha)^2 \S_{m\times n}(r^*)^2\ \p\pa{\{\hat r\geq r^*\}\cap \{\kf{E}{r^*}\geq (1-\alpha)\S_{m\times n}(r^*)\}}.
\end{eqnarray*}
Combining the analysis above with Gaussian concentration inequality for $\kf{E}{r^*}$, we have
$$\p\pa{\{\hat r\geq r^*\}\cap \{\kf{E}{r^*}\geq (1-\alpha)\S_{m\times n}(r^*)\}}\geq 1-2e^{\alpha^2/2}{e^{-\alpha^2\max(n,m)/2}\over 1-e^{-\alpha^2\max(n,m)/2}}.$$
We finally obtain the lower bound on the risk
$$\E\cro{\|X\widehat A-XA\|^2}\geq (1-\alpha)^2r^*(\max(n,m)-1)\pa{1-2e^{\alpha^2/2}{e^{-\alpha^2\max(n,m)/2}\over 1-e^{-\alpha^2\max(n,m)/2}}},$$
which is not compatible with the upper bound $c(K)$ that we would have if~(\ref{oracle1}) were also true with $K<1$.

When $q<m$, we start from $\|Y-X\widehat A_{r}\|^2=\|Y-PY\|^2+\|PY-X\widehat A_{r}\|^2$ with $P=X(X^*X)^+X^*$ and follow the same lines, replacing everywhere $E$ by $PE$ and $m$ by $q$.

\subsection{Proof of Proposition~\ref{propminimal2}}
As in the proof of Proposition~\ref{propminimal1}, we restrict for simplicity to the case where $\sigma^2=1$ and $q=m$, the general case being treated similarly. We write  $\penn(r)=\pen'(r)/(nm)$ and
for any integer $r^*\in [\min(n,m)/2,\min(n,m)-1]$, we set
$$\Omega_{*}=\ac{\kf{E}{r^*}\geq (1-\alpha) \E[\kf{E}{r^*}]}\bigcap_{r=1}^{\min(n,m)} \ac{\kf{E}{r}\leq (1+\alpha)\E[\kf{E}{r}]}.$$
According to the Gaussian concentration inequality we have
\begin{eqnarray*}
\p(\Omega_{*})&\geq& 1-2\sum_{r=1}^{\min(n,m)}e^{-\alpha^2\E[\kf{E}{r}]^2/2}\\
&\geq& 1-2e^{\alpha^2/2}\sum_{r=1}^{\min(n,m)}e^{-\alpha^2r\max(n,m)/2}
\end{eqnarray*}
where the last bound follows from Lemma~\ref{thm singular}.
For any $r\leq r^*$, we have on $\Omega_{*}$
\begin{eqnarray*}
\crit'(r^*)-\crit'(r)&=&\|E\|^2(\penn(r^*)-\penn(r))+\kf{E}{r}^2(1+\penn(r))-\kf{E}{r^*}^2(1+\penn(r^*))\\
&\leq& (1+\alpha)^2(\E[\|E\|]^2(\penn(r^*)-\penn(r))+\E[\kf{E}{r}]^2(1+\penn(r))(1+\alpha)^2\\
&&-\E[\kf{E}{r^*}]^2(1+\penn(r^*))(1-\alpha)^2.
\end{eqnarray*}
Since $\E[\|E\|]^2\leq nm=K\E[\kf{E}{r}]^2(1+\penn(r))/\penn(r)$, we have
\begin{eqnarray*}
\crit'(r^*)-\crit'(r)&\leq& (1+\alpha)^2(1-K)(1+\penn(r))\E[\kf{E}{r}]^2\\
&&-((1-\alpha)^2-(1+\alpha)^2K)(1+\penn(r^*))\E[\kf{E}{r^*}]^2\\
&\leq& (1+\alpha)^2(1-K)(1+\penn(r^*))\cro{\E[\kf{E}{r}]^2-(1-(1+\alpha)^{-2})\E[\kf{E}{r^*}]^2}.
\end{eqnarray*}
To conclude, we note that $\E[\kf{E}{r}]^2<r(\sqrt{n}+\sqrt{m})^2$, $\E[\kf{E}{r^*}]^2\geq (nm-1)/2$ and  $1-(1+\alpha)^{-2}\geq \alpha$, so the term in the bracket is smaller than
$$r(\sqrt{n}+\sqrt{m})^2-{1-K\over 8}(nm-1)$$
which is negative when $r\leq {1-K\over 8}\times{nm-1\over (\sqrt{n}+\sqrt{m})^2}$.

\subsection{Minimax rate : proof of Fact 1}\label{proof-minimax}

Let $X=U\Sigma V^*$ be a SVD decomposition of $X$, with the diagonal elements of $\Sigma$ ranked in decreasing order. Write $U_{q}$ and $V_{q}$ for the matrices derived from $U$ and $V$ by keeping the $q$-first columns, and $\Sigma_{q}$ for $q\times q$ upper-left block of $\Sigma$ (with notations as in {\tt R}, 
  $U_{q}=U[\ ,1:q]$, $V_{q}=V[\ ,1:q]$ and $\Sigma_{q}=\Sigma[1:q,1:q]$). We have $X=U_{q}\Sigma_{q} V^*_{q}$ and 
  $$Y=ZB+\sigma E,\quad \textrm{with}\ Z=U_{q}\Sigma_{q}\in {\bf R}^{m\times q}\ \textrm{and}\ B=V_{q}^*A \in {\bf R}^{q\times n}.$$
  Let $\tilde A$ be an arbitrary estimator of $A$ and set $\tilde B=V_{q}^*\tilde A$.
  Write $Z_{i}^*$ for the $i$th row of $Z$ and $\{e_{1},\ldots,e_{n}\}$ for the canonical basis of ${\bf R}^n$.
According to~(\ref{condminimax}), the map
$$B \to \mathcal{L}(B)=\cro{<Z_{i}e_{a}^*,B>/\sqrt{mn}\,}_{\begin{scriptsize}\begin{array}{l}i=1,\ldots,m\\ a=1,\ldots,n\end{array}\end{scriptsize}}$$
fulfills the Restricted Isometry condition RI$(r,\nu)$ of Rohde and Tsybakov~\cite{RohdeTsybakov} for all $r\leq \min(n,q)$ with
$$\nu^2={2mn\over \sigma_{1}(X)^2+\sigma_{q}(X)^2}\quad\textrm{and}\quad \delta={1-\rho^2\over 1+\rho^2}<1.$$
Theorem 2.5 in \cite{RohdeTsybakov} (with $\alpha=1/10$ and $\Delta=+\infty$) then ensures that there exists some constant $c_{\rho}>0$ depending only on $\rho$ such that
$$\inf_{\tilde B}\sup_{B\,:\,\mathrm{rank}(B)\leq r}\E[\|Z\tilde B-ZB\|^2]\geq 2c_{\rho}(q+n)r\sigma^2,\quad\textrm{for all }r\leq\min(n,q).$$
Let $B'$ be such that $\E[\|Z\tilde B-ZB'\|^2]\geq c_{\rho}(q+n)r\sigma^2$ and rank$(B')\leq r$. The matrix $A'=V_{q}B'\in{\bf R}^{p\times n}$ fulfills rank$(A')\leq r$ and
$$\E\cro{\|X\tilde A-XA'\|^2} = \E\cro{\|Z\tilde B-ZB'\|^2}\geq c_{\rho}(q+n)r\sigma^2.$$
In conclusion, for any $X$ fulfilling~(\ref{condminimax}), any estimator $\tilde A$ and any $r\leq \min(n,q)$, we have
$$\sup_{A\,:\, \mathrm{rank}(A)\leq r} \E\cro{\|X\tilde A-XA\|^2} \geq c_{\rho}(q+n)r\sigma^2.$$

\appendix
\section{Monte Carlo evaluation of $\S_{q\times n}(r)$}
\begin{verbatim}
SMonteCarlo <- function(q,n,Nsim=200){
   Sk <- array(0,c(Nsim,min(q,n)))	
   for (is in 1:Nsim) {
      s <- svd(matrix(rnorm(q*n),nrow=q,ncol=n),nu=0,nv=0)$d
      Sk[is,]<-sqrt(cumsum(s**2))
   }
   return(apply(Sk,2,mean))		
}
\end{verbatim}

\section{Marchenko-Pastur approximation of $\S_{q\times n}(r)$}
\begin{verbatim}
SMarchenkoPastur <-function(q,n,eps=10**(-9)){
   beta <- min(n,q)/max(n,q)
   alpha <- (1:min(n,q))/min(n,q)
   s<-rep(0,min(n,q))
   f <- function(x){
      return(sqrt((x-(1-sqrt(beta))^2)*((1+sqrt(beta))^2-x))/(2*pi*beta*x))}
   xf <- function(x){
      return(sqrt((x-(1-sqrt(beta))^2)*((1+sqrt(beta))^2-x))/(2*pi*beta))}
   for (a in 1:length(alpha)){
      m <- (1-sqrt(beta))^2
      M <- (1+sqrt(beta))^2
      while ((M-m)>eps) {
         if (integrate(f,(m+M)/2,(1+sqrt(beta))^2)$value<alpha[a]) M<-(m+M)/2
         else m<-(m+M)/2
      }
      s[a] <- integrate(xf,(m+M)/2,(1+sqrt(beta))^2)$value	
   }
   return(sqrt(s*n*q))
} 
\end{verbatim}

\end{document}